\documentstyle[amsfonts,bezier]{article}

\makeatletter
\def\eqnarray{%
  \stepcounter{equation}%
  \let\@currentlabel=\theequation
  \global\@eqnswtrue
  \global\@eqcnt\z@
  \tabskip\@centering
  \let\\=\@eqncr
  $$\halign to \displaywidth\bgroup\@eqnsel\hskip\@centering
  $\displaystyle\tabskip\z@{##}$&\global\@eqcnt\@ne
  \hfil$\displaystyle{{}##{}}$\hfil
  &\global\@eqcnt\tw@$\displaystyle\tabskip\z@{##}$\hfil
  \tabskip\@centering&\llap{##}\tabskip\z@\cr}
\makeatother

\makeatletter
  \renewcommand{\theequation}{%
        \thesection.\arabic{equation}}
  \@addtoreset{equation}{section}
\makeatother

\begin{document}

\newtheorem{th}{Donotwrite}[section]

\newtheorem{theorem}[th]{Theorem}
\newtheorem{proposition}[th]{Proposition}
\newtheorem{conjecture}[th]{Conjecture}
\newtheorem{lemma}[th]{Lemma}
\newtheorem{corollary}[th]{Corollary}
\newtheorem{remark}[th]{Remark}
\newtheorem{example}[th]{Example}

\newfont{\germ}{eufm10}
\newfont{\slsmall}{cmsl8}

\def\Ae{A^{(2)}_{2n}}
\def\Ao{A^{(2)}_{2n-1}}
\def\Aod{A^{(2)\dagger}_{2n-1}}
\def\bbar{\overline{b}}
\def\btilde{\tilde{b}}
\def\C{C^{(1)}_n}
\def\cd{\cdots}
\def\eps{\epsilon}
\def\et#1{\tilde{e}_{#1}}
\def\ft#1{\tilde{f}_{#1}}
\def\geh{\goth{g}}
\def\goth#1{\mbox{\germ #1}}
\def\id{\mbox{\sl id}\,}
\def\L{{\cal L}}
\def\La{\Lambda}
\def\la{\lambda}
\def\ol#1{\overline{#1}}
\def\ot{\otimes}
\def\P{{\cal P}}
\def\pbar{\overline{p}}
\def\Pcl{P_{cl}}
\def\Pcll{(P_{cl}^+)_l}
\def\Proof{\noindent{\sl Proof.}\quad}
\def\Q{{\bf Q}}
\def\q{{\bf q}}
\def\qed{~\rule{1mm}{2.5mm}}
\def\Rc{\check{R}}
\def\slchap{\widehat{\goth{sl}}_n}
\def\sln{\goth{sl}_n}
\def\Uq{U_q(\geh)}
\def\Uqp{U'_q(\geh)}
\def\veps{\varepsilon}
\def\vphi{\varphi}
\def\wt{\mbox{\sl wt}\,}
\def\wts{\mbox{\slsmall wt}\,}
\def\Z{{\bf Z}}
\def\Zn{\Z_{\ge0}}

\def\ep{\epsilon}

\title{ $q$-Wedge Modules for Quantized Enveloping \\
        Algebras of Classical Type }

\author{
Naihuan Jing\thanks{
Department of Mathematics, North Carolina State University,
Raleigh, NC 27695-8205, USA},
Kailash C. Misra\hspace{-0.3mm}$^*$
and Masato Okado\thanks{
Department of Informatics and Mathematical Science,
Graduate School of Engineering Science,
Osaka University, Toyonaka, Osaka 560-8531, Japan}
}

\date{}
\maketitle

\begin{abstract}
We use the fusion construction in the twisted quantum
affine algebras to obtain a unified method to deform
the wedge product for classical Lie algebras. As a byproduct
we uniformly 
realize all non-spin fundamental modules for
quantized enveloping algebras of classical types, and show
that they admit natural crystal bases as   
modules for the (derived) twisted quantum affine algebra. These crystal
bases are parametrized in terms of the $q$-wedge products.
\end{abstract}

\setcounter{section}{-1}

\section{Introduction}
Let us start with recalling a basic fact in representation theory
of finite-dimen\-sion\-al simple Lie algebras, which may be found in
an elementary textbook of Lie algebras. Let $V$ be the vector
representation of the simple Lie algebra $\sln$. The $k$-fold
exterior space $\bigwedge^k V$ ($1\le k<n$) turns out to be a new
representation of $\sln$, and 
\[
\bigwedge^k V\simeq V(\omega_k)\qquad\mbox{as $\sln$-modules}.
\]
Here $V(\omega_k)$ stands for the irreducible finite-dimen\-sion\-al
representation with the $k$-th fundamental weight $\omega_k$
as the highest weight.

Let us now consider the quantized enveloping algebra $U_q(\sln)$
in the sense of Drinfel'd and Jimbo. We also have the fundamental
representation $V(\omega_k)$ for $U_q(\sln)$. We would like to
realize $V(\omega_k)$ as a {\em wedge}-like representation.
Since the coproduct of the Hopf algebra $U_q(\sln)$ is 
non-cocommutative, the construction becomes a bit tricky. Let
$V=\bigoplus_{i=1}^n\Q(q)v_i$ be the vector representation of
$U_q(\sln)$. Our problem can be stated as follows.
\begin{description}
\item[]Find a suitable vector subspace $W$ of $V\ot V$ such that
\item[](1) $W$ is invariant under the action of $U_q(\sln)$,
\item[](2) at $q=1$, $W$ is generated by the vectors $v_i\ot v_i$
           and $v_i\ot v_j+v_j\ot v_i\,(i\neq j)$,
\item[](3) ${\displaystyle V^{\ot k}/\sum_{i=0}^{k-2}V^i\ot W\ot
           V^{\ot(k-2-i)}\simeq V(\omega_k)}$ as $U_q(\sln)$-modules.
\end{description}
In the case of $U_q(\sln)$, the answer is given by
\[
W=\mbox{ span}\{v_i\ot v_i,v_i\ot v_j+qv_j\ot v_i(i>j)\}.
\]
An important remark here is that $W$ is identified with the image
of the so-called $R$-matrix $\Rc(z)$ acting on $V\ot V$ at $z=q^2$.
\begin{eqnarray*}
\Rc(z)&=&(1-q^2z)\sum_iE_{ii}\otimes E_{ii}
         +q(1-z)\sum_{i\ne j}E_{ij}\otimes E_{ji}\\
      &&+(1-q^2)\Bigl(\sum_{i>j}+z\sum_{i<j}\Bigr)E_{ii}\otimes E_{jj}.
\end{eqnarray*}
Here $E_{ij}$ denotes the matrix unit, {\it i.e.} $E_{ij}v_l
=\delta_{jl}v_i$.

The purpose of this article is to quantize the wedge module
$\bigwedge^kV$ of the vector representation $V$ for the simple Lie
algebras of types $B_n,C_n$ and $D_n$ \`a la fusion construction
(e.g. \cite{KMN}). A weird thing here is that one needs to use the
$R$-matrices for affine Lie algebras of twisted types $\Ae$ and $\Ao$.
Furthermore, one has to use different $R$-matrices of $\Ao$
corresponding to different choices of the 0-vertex
(See Figure 1). If we compliantly choose
the $R$-matrices for $B^{(1)}_n,\C$ or $D^{(1)}_n$, we cannot
obtain a $q$-analogue of the wedge modules, as seen in \cite{DO}
or the $\C$ case in this article. We have nevertheless included
this $\C$ case since it can be treated in a similar fashion.
All the modules constructed here are not only $U_q(\ol{\geh})$-modules 
but also $\Uqp$-modules (See Table 1), and are shown to have a crystal base
in the sense of Kashiwara \cite{K1}. Moreover, except for the case of $\C$
they appear to be perfect of level one, which will be discussed elsewhere. 

A byproduct of our construction is that we also have
a unified way to realize all non-spin fundamental modules
for quantum enveloping algebras of classical types.

\section{Preliminaries}
In this section we shall fix conventions for the quantized enveloping
algebra $\Uq$ and recall the expressions of $R$--matrices.

\subsection{Dynkin datum}
We follow the notations of the affine Lie algebra in \cite{Kac}.
Let $\La_i$ ($i=0,\cd,n$) be the fundamental weights. Let $\alpha_i$ and
$h_i=\alpha^\vee_i$ ($i=0,\cd,n$) be the simple roots and coroots.
The generator of the null roots is denoted by $\delta=
\sum_{i=0}^{n}a_i\alpha_i$ and the canonical central element by
$c=\sum_{i=0}^{n}a^\vee_ih_i$. $d$ stands for the degree operator.

We consider the affine Lie algebras $\geh=\Ae(n\ge2),\Ao(n\ge3),\C(n\ge2)$.
For the Dynkin diagram of $\Ao$, we distinguish two different labelings
of the vertices by $\Ao$ and $\Aod$. The Dynkin diagrams for these
affine Lie algebras are given in Figure 1.
\begin{figure}[ht]
\setlength{\unitlength}{0.8pt}
\hspace{1cm}
\begin{picture}(300,60)(-50,0)
\put(-60,16){$A^{(2)}_{2n}$}
\put(5,20){\circle{10}}
\put(10,18){\line(1,0){36}}
\put(10,22){\line(1,0){36}}
\put(50,20){\line(-2,1){13}}
\put(50,20){\line(-2,-1){13}}
\put(55,20){\circle{10}}
\put(60,20){\line(1,0){40}}
\multiput(110,20)(6,0){10}{\circle*{1}}
\put(173,20){\line(1,0){40}}
\put(218,20){\circle{10}}
\put(223,18){\line(1,0){36}}
\put(223,22){\line(1,0){36}}
\put(263,20){\line(-2,1){13}}
\put(263,20){\line(-2,-1){13}}
\put(268,20){\circle{10}}
\put(0,3){\makebox(10,10){0}}
\put(50,3){\makebox(10,10){1}}
\put(213,3){\makebox(10,10){n-1}}
\put(263,3){\makebox(10,10){n}}
\end{picture}
\vskip0.7cm
\hspace{1cm}
\begin{picture}(300,60)(-50,0)
\put(-60,16){$A^{(2)}_{2n-1}$}
\put(5,20){\circle{10}}
\put(10,20){\line(1,0){40}}
\put(55,20){\circle{10}}
\put(55,25){\line(0,1){40}}
\put(55,70){\circle{10}}
\put(60,20){\line(1,0){40}}
\put(105,20){\circle{10}}
\put(110,20){\line(1,0){40}}
\multiput(160,20)(6,0){10}{\circle*{1}}
\put(223,20){\line(1,0){40}}
\put(268,20){\circle{10}}
\put(277,18){\line(1,0){36}}
\put(277,22){\line(1,0){36}}
\put(273,20){\line(2,1){13}}
\put(273,20){\line(2,-1){13}}
\put(318,20){\circle{10}}
\put(62,65){\makebox(10,10){0}}
\put(0,3){\makebox(10,10){1}}
\put(50,3){\makebox(10,10){2}}
\put(100,3){\makebox(10,10){3}}
\put(263,3){\makebox(10,10){n-1}}
\put(313,3){\makebox(10,10){n}}
\end{picture}
\vskip0.7cm
\hspace{1cm}
\begin{picture}(300,60)(-50,0)
\put(-60,16){$A^{(2)\dagger}_{2n-1}$}
\put(5,20){\circle{10}}
\put(10,18){\line(1,0){36}}
\put(10,22){\line(1,0){36}}
\put(50,20){\line(-2,1){13}}
\put(50,20){\line(-2,-1){13}}
\put(55,20){\circle{10}}
\put(60,20){\line(1,0){40}}
\put(105,20){\circle{10}}
\put(110,20){\line(1,0){40}}
\multiput(160,20)(6,0){10}{\circle*{1}}
\put(223,20){\line(1,0){40}}
\put(268,20){\circle{10}}
\put(268,25){\line(0,1){40}}
\put(268,70){\circle{10}}
\put(273,20){\line(1,0){40}}
\put(318,20){\circle{10}}
\put(0,3){\makebox(10,10){0}}
\put(50,3){\makebox(10,10){1}}
\put(100,3){\makebox(10,10){2}}
\put(263,3){\makebox(10,10){n-2}}
\put(275,65){\makebox(10,10){n}}
\put(313,3){\makebox(10,10){n-1}}
\end{picture}
\vskip-0.3cm
\hspace{1cm}
\begin{picture}(300,60)(-50,0)
\put(-60,16){$C^{(1)}_n$}
\put(5,20){\circle{10}}
\put(10,18){\line(1,0){36}}
\put(10,22){\line(1,0){36}}
\put(50,20){\line(-2,1){13}}
\put(50,20){\line(-2,-1){13}}
\put(55,20){\circle{10}}
\put(60,20){\line(1,0){40}}
\multiput(110,20)(6,0){10}{\circle*{1}}
\put(173,20){\line(1,0){40}}
\put(218,20){\circle{10}}
\put(227,18){\line(1,0){36}}
\put(227,22){\line(1,0){36}}
\put(223,20){\line(2,1){13}}
\put(223,20){\line(2,-1){13}}
\put(268,20){\circle{10}}
\put(0,3){\makebox(10,10){0}}
\put(50,3){\makebox(10,10){1}}
\put(213,3){\makebox(10,10){n-1}}
\put(263,3){\makebox(10,10){n}}
\end{picture}
\caption{Dynkin diagrams}
\end{figure}
The marks $(a_i)_{0\le i\le n}$ and comarks $(a^\vee_i)_{0\le i\le n}$ are
given by
\[
\begin{array}{rlrll}
(a_i)&=(1,2,\cd,2,2),&(a^\vee_i)&=(2,2,\cd,2,1)&\mbox{ for }\geh=\Ae,\\
&=(1,1,2,2,\cd,2,1),&&=(1,1,2,2,\cd,2,2)&\mbox{ for }\geh=\Ao,\\
&=(1,2,2,\cd,2,1,1),&&=(2,2,2,\cd,2,1,1)&\mbox{ for }\geh=\Aod,\\
&=(1,2,\cd,2,1),&&=(1,1,\cd,1,1)&\mbox{ for }\geh=\C.
\end{array}
\]
It is convenient to introduce an orthogonal basis $\{\eps_j\}_{1\le j\le n}$
with which the simple roots are given as follows:
\begin{eqnarray*}
\alpha_0&=&\delta-2\eps_1\quad\mbox{for }\geh=\Ae,\Aod,\C,\\
&=&\delta-\eps_1-\eps_2\quad\mbox{for }\geh=\Ao,\\
\alpha_i&=&\eps_i-\eps_{i+1}\quad(1\le i\le n-1),\\
\alpha_n&=&\eps_n\quad\mbox{for }\geh=\Ae,\\
&=&2\eps_n\quad\mbox{for }\geh=\Ao,\C,\\
&=&\eps_{n-1}+\eps_n\quad\mbox{for }\geh=\Aod.
\end{eqnarray*}
The invariant bilinear form $(\,|\,)$ in \cite{Kac} is  normalized so that
$(\theta|\theta)=2a^\vee_0$,
where $\theta=\delta-\alpha_0$. (Note that $a_0=1$ in our cases.)
With this normalization we have $(\eps_i|\eps_j)={\frac12}\delta_{ij}$
for $\geh=\C$, $\delta_{ij}$ for the other cases.

\subsection{Quantized enveloping algebra}
We recall the definition of the quantized enveloping algebra $\Uq$. Set
$P^*=\bigoplus_{i=0}^n\Z h_i\oplus\Z d$. $\Uq$ associated to $\geh$ is
the associative algebra generated by the symbols $e_i,f_i$ ($i=0,\cd,n$)
and $q^h$ ($h\in P^*$) satisfying the following relations:
\begin{eqnarray*}
&&q^0=1,\quad q^hq^{h'}=q^{h+h'},\\
&&q^he_iq^{-h}=q^{\langle h,\alpha_i\rangle}e_i,\quad
q^hf_iq^{-h}=q^{-\langle h,\alpha_i\rangle}f_i,\\
&&[e_i,f_j]=\delta_{ij}\frac{t_i-t_i^{-1}}{q_i-q_i^{-1}}
\end{eqnarray*}
where $q_i=q^{(\alpha_i|\alpha_i)/2},t_i=q^{(\alpha_i|\alpha_i)h_i/2}$,
\[
\sum_{k=0}^b(-1)^ke_i^{(k)}e_je_i^{(b-k)}=
\sum_{k=0}^b(-1)^kf_i^{(k)}f_jf_i^{(b-k)}=0\quad(i\neq j)
\]
where $b=1-\langle h_i,\alpha_j\rangle$. Here we set
$[k]_i=\left(q_i^k-q_i^{-k}\right)/\left(q_i-q_i^{-1}\right),
[k]_i!=\prod_{m=1}^k[m]_i,e_i^{(k)}=e_i^k/[k]_i!,f_i^{(k)}=f_i^k/[k]_i!$.

There are several coproducts of $\Uq$. We use the `lower' one
$\Delta=\Delta_{-}$ given by
\begin{eqnarray}
\Delta(e_i)&=&e_i\ot t_i^{-1}+1\ot e_i, \label{eq:cop}\\
\Delta(f_i)&=&f_i\ot 1+t_i\ot f_i, \nonumber\\
\Delta(q^h)&=&q^h\ot q^h. \nonumber
\end{eqnarray}
For the relations among different coproducts, see \cite{KMPY} for example.
We also define $\Uqp$ as the subalgebra of $\Uq$ generated by the elements
$e_i,f_i,t_i$ ($i=0,\cdots,n$), and $U_q(\ol{\geh})$ as the subalgebra
generated by $e_i,f_i,t_i$ ($i=1,\cdots,n$). $U_q(\ol{\geh})$ is canonically
viewed as the quantized enveloping algebra associated to the simple Lie
algebra $\ol{\geh}$.
\begin{table}
\[
\begin{array}{cccccc}
\geh&=&\Ae&\Ao&\Aod&\C\\
\overline{\geh}&=&B_n& C_n& D_n& C_n
\end{array}
\]
\caption{Associated simple Lie algebras}
\end{table}

\subsection{Representations and the $R$--matrix}
Define an index set $J$ by
\begin{eqnarray*}
J&=&\{0,\pm1,\cdots,\pm n\}\quad\hbox{for }\geh=\Ae,\\
 &=&\{\pm1,\cdots,\pm n\}\quad\hbox{for }\geh=\Ao,\Aod,\C.
\end{eqnarray*}
We introduce a linear order $\prec$ in $J$ by
\[
1\prec2\prec\cdots\prec n\,(\prec0)\prec-n\prec\cdots\prec-2\prec-1.
\]
We now consider the `vector' representation
$(\pi,V)$ of $U'_q(\geh)$.
Let $\{v_j\mid j\in J\}$ be the lower global crystal base \cite{K2}.
Denoting the matrix units by $E_{ij}$ {\it i.e.} $E_{ij}v_k=\delta_{jk}v_i$,
the actions of the generators read as follows. ($\pi(f_i)=\pi(e_i)^t$,
unless otherwise stated.)
\begin{eqnarray*}
\pi(e_0)&=&E_{-1,1}\quad\hbox{for }\geh=\Ae,\Aod,\C,\\
        &=&E_{-1,2}+E_{-2,1}\quad\hbox{for }\geh=\Ao,\\
\pi(t_0)&=&\sum_{j\in J}q^{-2\delta_{j1}+2\delta_{j,-1}}E_{jj}
          \quad\hbox{for }\geh=\Ae,\Aod,\C,\\
        &=&\sum_{j\in
J}q^{-\delta_{j1}-\delta_{j2}+\delta_{j,-1}+\delta_{j,-2}}
          E_{jj}\quad\hbox{for }\geh=\Ao,\\
\pi(e_i)&=&E_{i,i+1}+E_{-i-1,-i}\quad(1\le i\le n-1),\\
\pi(t_i)&=&\sum_{j\in J}q^{\delta_{ji}-\delta_{j,i+1}+\delta_{j,-i-1}
          -\delta_{j,-i}}E_{jj}\quad(1\le i\le n-1),\\
\pi(e_n)&=&[2]_nE_{n0}+E_{0,-n}\quad\hbox{for }\geh=\Ae,\\
        &=&E_{n,-n}\quad\hbox{for }\geh=\Ao,\C,\\
        &=&E_{n-1,-n}+E_{n,-n+1}\quad\hbox{for }\geh=\Aod,\\
\pi(f_n)&=&E_{0n}+[2]_nE_{-n,0}\quad\hbox{for }\geh=\Ae,\\
\pi(t_n)&=&\sum_{j\in J}q^{\delta_{jn}-\delta_{j,-n}}E_{jj}
          \quad\hbox{for }\geh=\Ae,\\
        &=&\sum_{j\in J}q^{2\delta_{jn}-2\delta_{j,-n}}E_{jj}
          \quad\hbox{for }\geh=\Ao,\C,\\
        &=&\sum_{j\in
J}q^{\delta_{j,n-1}+\delta_{jn}-\delta_{j,-n}-\delta_{j,-n+1}}
          E_{jj}\quad\hbox{for }\geh=\Aod.
\end{eqnarray*}

Let us define the so-called `evaluation' representation $(\pi_z,V_z)$
associated with $(\pi,V)$  by setting
$V_z=V[z,z^{-1}],\pi_z(e_i)=z^{\delta_{i0}}\pi(e_i),
\pi_z(f_i)=z^{-\delta_{i0}}\pi(f_i),\pi_z(t_i)=\pi(t_i)$.
Let $\Rc(z)$ be the $R$--matrix, which satisfies
\begin{equation} \label{eq:int_rel}
\Rc(z_1/z_2)(\pi_{z_1}\otimes\pi_{z_2})\Delta(x)
=(\pi_{z_2}\otimes\pi_{z_1})\Delta(x)\Rc(z_1/z_2)
\quad\hbox{for all}~x\in\Uqp.
\end{equation}
The explicit expressions for our $\Rc(z)$ have been known in \cite{J1},
though we need a slight modification due to a different choice of the
basis and the coproduct, except for the $\Ao$ case. Note that the expressions
of $\Rc(z)$ for $\Ao$ and $\Aod$ should be different due to different
choices of the 0-vertex. The expression of $\Rc(z)$ 
for $\Ao$ seems to be new.
Up to an overall function of $z$, $\Rc(z)$ is uniquely determined by
\begin{eqnarray*}
\Rc(z)&=&(1-q^2z)(1-\xi z)\sum_{i\ne0}E_{ii}\otimes E_{ii}\\
     &&+q(1-z)(1-\xi z)\sum_{i\ne\pm j}E_{ij}\otimes E_{ji}\\
     &&+(1-q^2)(1-\xi z)\Bigl(\sum_{i\succ j,i\ne-j}+z\sum_{i\prec j,i\ne-j}
       \Bigr)E_{ii}\otimes E_{jj}\\
     &&+\sum_{i,j} a_{ij}(z)E_{-i,j}\otimes E_{i,-j}.
\end{eqnarray*}
Here
\begin{equation} \label{eq:def_a}
a_{ij}(z)=\left\{
\begin{array}{ll}
         (q^2-\xi z)(1-z)+\delta_{i0}(1-q)(q+z)(1-\xi z)&(i=j)\\
         (1-q^2)(\veps_j\veps_i^{-1}(-q)^{\overline{j}-\overline{i}}(z-1)
                 +\delta_{i,-j}(1-\xi z))
               &(i\prec j)\\
         (1-q^2)z(\veps_j\veps_i^{-1}\xi (-q)^{\overline{j}-\overline{i}}(z-1)
                  +\delta_{i,-j}(1-\xi z))
               &(i\succ j).
\end{array}\right.
\end{equation}
Here $\veps_j=1$ ($j\in J$) for $\geh=\Aod$ and $\veps_j=1$ ($j>0$),
$[2]_n$ ($j=0$), $-1$ ($j<0$) for the other cases.
Also $\overline{j}$ is defined by
\[
\overline{j}=\left\{
\begin{array}{ll}
        j&(j=1,\cdots,n)\\
        n+1&(j=0)\\
        j+N&(j=-n,\cdots,-1),
\end{array}\right.
\]
and $N,\xi$ are given below.
\[
\begin{array}{cccccc}
\geh&=&\Ae&\Ao&\Aod&\C\\
N&=&2n+1&2n+2&2n&2n\\
\xi&=&-q^{2n+1}&-q^{2n}&-q^{2n}&q^{2n+2}
\end{array}
\]


\section{$q$-wedge constructions}

In this section we will use the fusion construction to
deform the wedge product. The basis for the deformed wedge product
relies on the
$R$-matrix and its spectral decomposition. 

\subsection{$R$-matrices and their spectral decompositions}

Let $\ol{\la}$ be the classical part of a weight $\la$ \cite{Kac}.
Let $V(\ol{\la})$ be the highest weight $U_q(\ol{\geh})$-module
with highest weight $\ol{\la}$.
Recall $V(=V(\ol{\La}_1)$ as a $U_q(\ol{\geh})$-module)
is the vector representation of the quantized
enveloping algebra $\Uqp$. Direct computation gives the following proposition.

\begin{proposition} \label{P:tensordec}
As a $U_q(\ol{\geh})$-module, the tensor product $V^{\otimes 2}$ decomposes
itself as follows.
\[
V(\ol{\La}_1)\otimes V(\ol{\La}_1)=V(2\ol{\La}_1)\oplus V(\ol{\La}_2)
\oplus V(0),
\]
and the highest weight vectors are respectively:
\begin{eqnarray*}
u_{2\ol{\La}_1}&=&v_1\otimes v_1, \\
u_{\ol{\La}_2}&=&v_1\otimes v_2-q v_2\otimes v_1,          \\
u_{0}
&=&\sum_{i=1}^{n}(-q)^{i-1}v_i\otimes v_{-i}
-\sum_{i=1}^{n}(-q)^{-i-1}\xi' v_{-i}\otimes v_i
-\frac{(-q)^{n-1}}{[2]_n}v_0\otimes v_0,\\
\end{eqnarray*}
where we take $v_0=0$ except for type $\Ae$ and $\xi'=\xi$ except
$\xi'=q^{2n+2}$ for type $\Ao$.
\end{proposition}

\begin{proposition} \label{P:spectral}
Let $P_{2\ol{\La}_1}$, $P_{\La_2}$ and $P_{0}$ be the projections from $V(\ol{\La}_1)\otimes
V(\ol{\La}_1)$
to the corresponding $U_q(\overline{\goth g})$-irreducible components in 
Prop. (\ref{P:tensordec}) respectively. Then we have the following spectral decomposition
of $\Rc(z)$.
\[
\Rc(z)=(1-q^2z)(1-\xi z)P_{2\ol{\La}_1}\oplus (1-\xi z)(z-q^2) P_{\ol{\La}_2}
\oplus (1-q^2z)(z-\xi)P_0
\]
\begin{flushright}
for $\geh=\Ae, \Aod, \C$,\phantom{flush}
\end{flushright}
and
\[
\Rc(z)=(1-q^2z)(1-\xi z)P_{2\ol{\La}_1}\oplus (1-\xi z)(z-q^2) P_{\ol{\La}_2}
\oplus (z-q^2)(z-\xi)P_0
\]
\begin{flushright}
for $\geh=\Ao$.\phantom{flush}
\end{flushright}
\end{proposition}

\noindent
{\it Proof.} Since $\Rc(z)$ commutes with the action of $U_q(\ol{\geh})$,
it is a constant on each irreducible component.
{}From Proposition \ref{P:tensordec} the
$R$-matrix $\Rc(z)$ can be written as:
$$
\Rc(z)=c_{2\ol{\La}_1}P_{2\ol{\La}_1}+c_{\ol{\La}_2}P_{\ol{\La}_2}+c_0P_0.
$$
Since the other cases are found in the literature (see \cite{KMN} for
$\geh=\Ae,\Ao$ and \cite{J2} for $\geh=\C$), we prove this only in the
cases of $\Aod$ and $\Ao$. One computes
\begin{eqnarray*}
&&(\pi_{z_1}\ot\pi_{z_2})
\Delta(f_0f_1\cdots f_{n-2}f_nf_{n-1}\cdots f_2)u_{\ol{\La}_2}=
(q^{-1}z_2^{-1}-qz_1^{-1})u_{2\ol{\La}_1},\\
&&(\pi_{z_1}\ot\pi_{z_2})
\Delta(f_0)u_0=(q^{-2}z_2^{-1}+q^{2n-2}z_1^{-1})u_{2\ol{\La}_1}, \qquad
\mbox{for $\Aod$}; and \\
&&(\pi_{z_1}\ot\pi_{z_2})
\Delta(f_0f_2\cdots f_{n-1}f_nf_{n-1}\cdots f_2)u_{\ol{\La}_2}=
(q^{-1}z_2^{-1}-qz_1^{-1})u_{2\ol{\La}_1},\\
&&(\pi_{z_1}\ot\pi_{z_2})
\Delta(f_0)u_0=(q^{-1}z_2^{-1}+q^{2n-1}z_1^{-1})u_{2\ol{\La}_1}, \qquad
\mbox{for $\Ao$}.
\end{eqnarray*}
Using the intertwining relation (\ref{eq:int_rel}) we get
\begin{eqnarray*}
c_{2\ol{\La}_1}(q^{-1}-qz^{-1})&=&c_{\ol{\La}_2}(q^{-1}z^{-1}-q),\\
c_{2\ol{\La}_1}(q^{-2}+q^{2n-2}z^{-1})&=& c_0(q^{-2}z^{-1}+q^{2n-2}),
\qquad\mbox{for $\Aod$};\\
c_{2\ol{\La}_1}(q^{-1}-qz^{-1})&=&c_{\ol{\La}_2}(q^{-1}z^{-1}-q),\\
c_{2\ol{\La}_1}(q^{-1}+q^{2n-1}z^{-1})&=& c_0(q^{-1}z^{-1}+q^{2n-1}),
\qquad\mbox{for $\Aod$}
\end{eqnarray*}
where $z=z_1/z_2$. These relations
determine the ratio of the coefficients $c_{2\ol{\La}_1}:c_{\ol{\La}_2}:c_0$
uniquely.
\hspace{\fill} $\Box$

\subsection{The $q$-wedge relations}

Let $W$ be the image of $\Rc(q^2)$ in $V^{\ot2}$, then
{}from the spectral decomposition of the $R$-matrix
(Proposition \ref{P:spectral}) it easily follows that
$$
W=Im\,\Rc(q^2)=Ker\,\Rc(q^{-2})\simeq V^{\ot2}/ Ker\,\Rc(q^2).
$$

\begin{proposition} \label{P:relations}
$W$ is generated by the vectors
\begin{eqnarray} \label{E:rel-1}
&&v_i\ot v_i~(i\neq0),\quad v_i\ot v_j+qv_j\ot v_i~(i\succ j,i\neq\pm j),\\ \label{E:rel-2}
&&v_{-i}\ot v_i+q^2v_i\ot v_{-i}+q(v_{i+1}\ot v_{-i-1}+v_{-i-1}\ot v_{i+1})~
(1\le i<n),
\end{eqnarray}
plus the following additional ones:
\begin{eqnarray*}
&&v_{-1}\ot v_1+v_1\ot v_{-1},\quad
v_{-n}\ot v_n+q^2v_n\ot v_{-n}+q^{\frac12}v_0\ot v_0\quad
\mbox{for }\geh=\Ae,\\ 
&&v_{-n}\ot v_n+q^2v_n\ot v_{-n}\quad
\mbox{for }\geh=\Ao,\\ 
&&v_{-1}\ot v_1+v_1\ot v_{-1}\quad
\mbox{for }\geh=\Aod,\\ 
&&v_{-1}\ot v_1+v_1\ot v_{-1},\quad
v_{-n}\ot v_n+q^2v_n\ot v_{-n}\quad
\mbox{for }\geh=\C.
\end{eqnarray*}
\end{proposition}
{\it Proof.}
The $R$ matrix $\Rc(z)$ preserves the weight, thus it is enough to
show that the listed vectors generate the spanning set
$\{ \Rc(q^2) (v_i\otimes v_j) \mid i, j\in J\}$.

For $i\in J\setminus\{0\}$ we have
$$
\Rc(z)(v_i\otimes v_i)=(1-q^2z)(1-\xi z)v_i\otimes v_i,
$$
which implies that $v_i\otimes v_i\in W$.

For $i\neq\pm j$ we have
\begin{eqnarray*}
\Rc(q^2)(v_i\otimes v_j)&=&q(1-q^2)(1-\xi q^2)v_j\otimes v_i\\
&&+(1-\xi q^2)(1-q^2)\times
\left\{
\begin{array}{rl}
v_i\ot v_j&\mbox{ if }i\succ j\\
q^2v_i\ot v_j&\mbox{ if }i\prec j
\end{array}\right.\\
&=& \left\{
\begin{array}{rl}
(1-q^2)(1-\xi q^2)(q v_j\otimes v_i+v_i\otimes v_j) &\mbox{ if $i\succ j$}\\
q(1-q^2)(1-\xi q^2)(v_j\otimes v_i+qv_i\otimes v_j) &\mbox{ if $i\prec j$},
\end{array}\right.
\end{eqnarray*}
thus $v_i\otimes v_j+q v_j\otimes v_i (i\succ j,  i\neq -j)\in W$.

As for the weight $0$ vectors $\Rc(q^2)(v_j\otimes v_{-j})$, we only
show it for the case of $\Aod$ for the purpose of
simplicity, since all other cases are similar.

Let $A_j=A_j(z)=\sum_{i\in J}a_{ij}(z)v_{-i}\otimes v_i=\Rc(z)(v_j\otimes v_{-j})$,
where $a_{ij}(z)$ are defined in (\ref{eq:def_a}).
Write $u_i=v_i\otimes v_{-i}$. We have for $1\leq j\leq n-1$
\begin{eqnarray*}
A_j+q^{-1}A_{j+1}&=&(z-1)q(z\xi-1)u_{-j-1}+(1-z)(1-\xi z)u_{-j}\\
&&\qquad +z(1-q^2)(1-\xi z)u_j+q^{-1}z(1-q^2)(1-\xi z)u_{j+1}\\
&=& (1-\xi z)\{(1-z)(qu_{-j-1}+u_{-j})+z(1-q^2)(u_j+q^{-1}u_{j+1})\}
\end{eqnarray*}
At $z=q^2$ we see that $u_{-j}+q^2u_j+q(u_{j+1}+u_{-j-1})\in W$.
Similarly we have at $z=q^2$
$$A_{-j-1}+q^{-1}A_{-j}=(1-q^2)(1-q^2\xi)\{q^{-1}u_{-j}+u_{-j-1}+u_{j+1}
+qu_{j}\},$$
which does not produce any new vectors in $W$.
Also we have
$$A_{n}-A_{-n}=(z-q^2)(1-\xi z)(u_n-u_{-n})$$
which vanishes when $z=q^2$.
These relations show that $W$ is generated by the vectors
$$
v_i\otimes v_i, v_i\otimes v_j +q v_j\otimes v_i (i \succ j, i\neq -j), \\
u_{i+1}+u_{-i-1}+(q^{-1}u_{-i}+qu_i),$$
plus one of the vectors $A_j(q^2)$.

We claim that for positive $j\in J$
\begin{eqnarray*}
A_j(q^2)&\equiv&q^{-2(n-j)}A_{-j}(q^2)\\
&\equiv&(-q)^{j-1}q^2(1-q^4)(u_1+u_{-1}) \quad
\mbox{mod} \,\, u_{i+1}+u_{-i-1}+(q^{-1}u_{-i}+qu_i).
\end{eqnarray*}
In fact we have
\begin{eqnarray*}
A_1&=& (q^2-\xi z)(1-z)u_{-1}
+(1-q^2)z\left(\xi(-q)^{2-2n}(z-1)+(1-\xi z)\right) u_1 \\
&&+\sum_{i=2}^n(1-q^2)z\xi(z-1)\left( (-q)^{1-i}u_{-i}+(-q)^{1-2n+i}u_{i}
\right) \\
&\equiv& (q^2-\xi z)(1-z)u_{-1}
+(1-q^2)z\left(\xi(-q)^{2-2n}(z-1)+(1-\xi z)\right) u_1 \\
&&+(1-q^2)z\xi(z-1)\left( (-q)^{-1}+(-q)^{-3}+\cd+(-q)^{3-2n}\right)
(-q^{-1}u_{-1}-qu_1) \\
&=&z(1-z^2)u_{-1}+z\left( (1-q^2z)+\xi(q^2-z)\right)u_1,
\end{eqnarray*}
which is equivalent to $q^2(1-q^4)(u_1+u_{-1})$ at $z=q^2$.
\hspace{\fill} $\Box$

Define the action of $x=e_i,f_i,t_i\in\Uqp$ on $V^{\ot k}$ by
\[
x(u_1\ot\cd\ot u_k)=(\pi_{z_1}\ot\cd\ot\pi_{z_k})\Delta^{(k)}(x)
(u_1\ot\cd\ot u_k).
\]
Here $z_j=(-q)^{2j-k}$ and 
$\Delta^{(k)}=(\Delta\ot\underbrace{\id\ot\cd\ot\id}_{k-2})\circ\cd\circ
(\Delta\ot\id)\circ\Delta$. Explicitly, it reads as 
\begin{eqnarray*}
e_i(u_1\ot u_2\ot \cdots\ot u_k)
&=&\sum_{j=1}^k u_1\ot \cdots \ot u_{j-1}
\ot z_j^{\delta_{i0}}e_iu_j\ot t_i^{-1}u_{j+1}\ot \cdots \ot
t_i^{-1}u_k\\
f_i(u_1\ot u_2\ot \cdots
\ot u_k)&=&\sum_{j=1}^k t_iu_1\ot \cdots \ot t_iu_{j-1}
\ot z_j^{-\delta_{i0}}f_iu_j\ot u_{j+1}\ot \cdots \ot u_k\\
t_i(u_1\ot u_2\ot \cdots
\ot u_k)&=&t_iu_1\ot \cdots \ot
t_iu_k.
\end{eqnarray*}
Because of the intertwining property (\ref{eq:int_rel}), it is immediate 
to see that the subspace $\sum_{i=1}^{k-1} V^{\otimes (i-1)}\otimes W\otimes
V^{\otimes(k-i-1)}$ is a $\Uqp$-submodule of $V^{\ot k}$. Thus setting 
$$
V^k=V^{\otimes k}/\sum_{i=1}^{k-1} V^{\otimes (i-1)}\otimes W\otimes
V^{\otimes(k-i-1)},
$$
$V^k$ turns out to be a $\Uqp$-module. 
For simplicity we will write the image of
$u_1\otimes u_2\otimes\cdots
\otimes u_k$ as $u_1\wedge u_2\wedge\cdots\wedge u_k$.

At $q=1$ the module $V^k$ degenerates into the wedge product
for the simple Lie algebra except for the cases when
the algebra is $C_n^{(1)}$ or $\Aod$ and $k=n$ . We will include some examples
at the end of this section.

We set $\omega_i=\epsilon_1+\cdots+\epsilon_i$ and
$\overline{\omega}_n=\epsilon_1+\cdots+\epsilon_{n-1}-
\epsilon_n$. The relations between $\omega_i$ and $\Lambda_i$ are
as follows.
\begin{eqnarray*}
\omega_i&=&\left\{\begin{array}{ll}
 2\Lambda_n, & i=n; B_n, D_n, \\
 \Lambda_{n-1}+\Lambda_n, & i=n-1; D_n,\\
\Lambda_i, & \mbox {otherwise}\\
\end{array}\right.\\
\overline{\omega}_{n}&=&\quad 2\Lambda_{n-1},  \qquad\qquad D_n.
\end{eqnarray*}

\begin{theorem}
For each $k$ the module $V^k$ is isomorphic to the highest weight
$U_q(\overline{\goth g})$-module
with the highest weight $\omega_k$ except for
\begin{eqnarray*}
V^{n}&=& V(\omega_{n})\oplus V(\overline{\omega}_n), \qquad
\mbox{for $\Aod$}, \\
V^{k}&=& V(\omega_{k})\oplus V(\omega_{k-2})\oplus\cdots
\oplus V(\omega_{k\,mod\,2}) , \qquad
\mbox{for $\Ao$}.
\end{eqnarray*}
\end{theorem}
{\it Proof.}
By the relations of $V^k$ we see that the space is generated by the
set of vectors:
$$v_{i_1}\wedge v_{i_2} \wedge \cdots \wedge v_{i_k},$$
where $i_1 \prec i_2 \prec \cdots \prec i_k$ are $k$ indices from the
set $J$. Therefore the dimension of the space is less than or
equal to $$\left(\begin{array}{c} |J|\\k \end{array}\right).$$

On the other hand the module $V^k$ contains the highest weight vector
$v_1 \wedge v_2 \wedge \cdots \wedge v_k$ of weight
$$\epsilon_1+\epsilon_2+\cdots +\epsilon_k=\omega_k.$$

When $k=n$ in the case of $\Aod$, there is another highest weight
vector $v_1 \wedge v_2 \wedge \cdots \wedge v_{n-1}\wedge v_{-n}$,
which is of the weight
$$\epsilon_1+\epsilon_2+\cdots +\epsilon_{n-1}-\epsilon_n=
\overline{\omega}_n.$$

When $k\le n$ in the case of $\Ao$, there are $[k/2]+1$ highest weight vectors.
For a multi-index set $I=(i_1, i_2, \cdots, i_l)$ denote
$$|I|=i_1+i_2+\cdots+ i_l.$$

We claim that in the case of $\Ao$, 
for each $l=0, 1, \cdots, [k/2]$ the following vector
$w_{k-2l}$ is a highest weight vector of weight
$\omega_{k-2l}=\ep_1+\cdots+\ep_{k-2l}$:
\begin{equation} \label{eq:hwv}
w_{k-2l}=\sum_{k-2l<i_1<i_2<\cdots <i_l}
(-q)^{|I|}v_1\wedge v_2\wedge\cdots\wedge v_{k-2l}\wedge v_{i_1}\wedge
\cdots\wedge v_{i_l}\wedge v_{-i_l}\wedge\cdots\wedge v_{-i_1}.
\end{equation}

In fact we can check directly that the action of
$e_i|_{V^k}=\sum 1\otimes \cdots 1\otimes e_i\otimes t_i^{-1}\otimes\cdots
\otimes t_i^{-1}$ ($i\neq0$)  kills the vector.
Since $wt(v_{\pm i})=\pm\epsilon_i$,
we have $wt(w_{k-2l})=\omega_{k-2l}$ .

First we consider the cases other than $C_n^{(1)}$.
By the dimension formula for irreducible modules of simple
Lie algebras \cite{OV} it follows that

\begin{eqnarray*}
\dim\,V^k&\geq & \dim\,V(\omega_k)=\left(\begin{array}{c} |J|\\k
\end{array}\right), \qquad\mbox{for  $\Aod (k\neq n), \Ae$,}\\
\dim\, V^k &\geq & \dim\, V(\omega_n)+\dim\,
V(\overline{\omega}_n)=\left(\begin{array}{c} |J|\\k \end{array}\right),
\qquad\mbox{for $\Aod (k=n)$},\\
\dim\, V^k&\geq & \dim\, V(\omega_k)+\dim\, V(\omega_{k-2})+\cdots +\dim\,
V(\omega_{k\,mod\,2})\\
&=&\left(\begin{array}{c} |J|\\k \end{array}\right),
\qquad\mbox{for $\Ao$}.
\end{eqnarray*}

In the $C_n^{(1)}$ case we need the formula
\begin{equation} \label{eq:formula}
\sum_{j=1}^n (-q)^j v_j\wedge v_{-j}=0,
\end{equation}
which is easily shown from the relations in Proposition \ref{P:relations} (cf. Prop. \ref{P:relations-1}).
Since the action of $e_i$ ($i\neq0$) is the same as in the $\Ao$ case,
(\ref{eq:hwv}) give rise to highest weight vectors if 
they are not zero. Moreover, 
from the properties of weights, they are the only possible highest weight vectors.
Thus it suffices to check $w_{k-2l}=0$ if $l>0$. Interchanging $i_l$ and
$i_{l-1}$ and using relations, we have 
\[
q^2w_{k-2l}=\sum_{k-2l<i_1<\cd<i_{l-2}<i_l<i_{l-1}}
(-q)^{|I|}v_{1,\cd,k-2l,i_1,\cd,i_l,-i_l,\cd,-i_1},
\]
where $v_{j_1,\cd,j_k}$ denotes $v_{j_1}\wedge\cd\wedge v_{j_k}$.
Similarly we have 
\[
q^{2l-2j}w_{k-2l}=\sum_{k-2l<i_1<\cd<i_{j-1}<i_l<i_j<\cd<i_{l-1}}
(-q)^{|I|}v_{1,\cd,k-2l,i_1,\cd,i_l,-i_l,\cd,-i_1}.
\]
If $i_l$ is equal to some other $i_j$, we see $v_{i_1,\cd,i_l}=0$. Thus
using (\ref{eq:formula}) we have
\begin{eqnarray*}
(1+q^2+\cd+q^{2l-2})w_{k-2l}\hspace{-1.4cm}&&\\
&=&\sum_{k-2l<i_1<\cd<i_{l-1}}\sum_{k-2l<i_l}
(-q)^{|I|}v_{1,\cd,k-2l,i_1,\cd,i_l,-i_l,\cd,-i_1}\\
&=&-\sum_{k-2l<i_1<\cd<i_{l-1}}\sum_{i_l\le k-2l}
(-q)^{|I|}v_{1,\cd,k-2l,i_1,\cd,i_l,-i_l,\cd,-i_1},
\end{eqnarray*}
which is zero unless $q^{2l}=1$.

Thus the theorem is proved.
\hspace{\fill} $\Box$

\subsection{Remarks and Examples}

We will use the $q$-wedge product to calculate
some analog of classical identities.
First we make several observations about the $q$-wedge products.

In the cases of $\Aod$, $\Ae$ and $C_n^{(1)}$, the vector
$$u_0=
\sum_{i=1}^{n}(-q)^{i-1}v_i\wedge v_{-i}
-\sum_{i=1}^{n}(-q)^{2n-i+1}v_{-i}\wedge v_i
-\frac{(-q)^{n-1}}{[2]_n}v_0\wedge v_0=0.
$$

In the case of $\Ao$, we have
\begin{eqnarray}\label{E:wt0}
u_0&=&\sum_{i=1}^{n}(-q)^{i-1}v_i\wedge v_{-i}
-\sum_{i=1}^{n}(-q)^{-i-1}\xi' v_{-i}\wedge v_i \nonumber\\
&=&(1+q^{2n+2})\sum_{i=1}^{n}(-q)^{i-1}v_i\wedge v_{-i}.
\end{eqnarray}
This is proved by an inductive calculation based on the wedge
relations. Let $v_0'$ be the vector $\sum_{i=1}^n
(-1)^{i-1}q^{-i+1}v_{-i}\wedge v_i$.

\begin{eqnarray*}
v_0'&=&v_{-1}\wedge v_1-q^{-1}v_{-2}\wedge v_2+\cdots +(-q)^{n-1}v_{-n}\wedge
v_n\\
&=&-q^2v_1\wedge v_{-1}-q(v_{2}\wedge v_{-2}+ v_{-2}\wedge v_{2})-q^{-2}
v_{-2}\wedge v_{2}+\cdots\\
&=&-q^2v_1\wedge v_{-1}+q^3v_2\wedge v_{-2}+q[2](v_{-3}\wedge v_{3}+
v_3\wedge v_{-3})+q^{-2}v_{-3}\wedge v_3+\cdots\\
&=&-q^{2}v_1\wedge v_{-1}+\cdots+(-1)^{i-1}q^iv_{i-1}\wedge v_{i-1}
+(-1)^{i-1}q[i-1](v_{-i}\wedge v_i\\
&&\qquad+v_{i}\wedge v_{-i})+(-q)^{-i+1}
v_{-i}\wedge v_i+\cdots+(-q)^{n-1}v_{-n}\wedge v_n
\end{eqnarray*}
where we assume inductively at the $i$-step. Using the relations
\begin{eqnarray*}
v_{-i}\wedge v_i&=&-q^2v_i\wedge v_{-i}-q(v_{-i-1}\wedge v_{i+1}+
v_{i+1}\wedge v_{-i-1}),\\
\ [i]&=&q[i-1]+q^{-i+1},
\end{eqnarray*}
we see by induction that
\begin{eqnarray*}
v_0'&=&-q^{2}v_1\wedge v_{-1}+\cdots+(-1)^{n-1}q^nv_{n-1}\wedge v_{n-1}
+(-1)^{n-1}q[n-1](v_{-n}\wedge v_n\\
&&\qquad +v_{n}\wedge v_{-n})+
\cdots+(-q)^{n-1}v_{-n}\wedge v_n\\
&=&-q^{2}v_1\wedge v_{-1}+\cdots+(-1)^{n-1}q^nv_{n-1}\wedge v_{n-1}\\
&&\qquad +(-1)^{n-1}(q[n-1]-q^2[n])v_n\wedge v_{-n}\\
&=&-q^2v_1\wedge v_{-1}+q^3 v_2\wedge v_{-2}+\cdots -(-q)^{n+1}v_n\wedge v_{-n}
\end{eqnarray*}
which proves (\ref{E:wt0}).

From the relations it is easy to obtain the following formula.
\[
v_{-i}\wedge v_i=-q^2v_i\wedge
v_{-i}+(1-q^2)\sum_{k=1}^{n-i}(-q)^kv_{i+k}\wedge v_{-(i+k)}
\]
for $i=1, \cdots, n-1$.

Note that except for the case of
$C_n^{(1)}$, the theorem implies that
the set of vectors $v_{i_1}\wedge v_{i_2}
\wedge \cdots \wedge v_{i_k}$, where $i_1 \prec i_2 \prec \cdots \prec i_k$
are $k$ indices from the set $J$, forms a basis for the module $V^k$. At $q=1$
the module $V^k$
specializes to the ordinary wedge representation of the subalgebra
$\overline{\geh}$ except for
$C_n^{(1)}$. In fact, it is
well-known that
the wedge representations for classical Lie algebras are as follows \cite{OV}:

\[
\bigwedge ^k V(\Lambda_1)\simeq
\left\{\begin{array}{lll} V(\omega_k) 
& &\mbox{ for $A_n, B_n, D_n$ and $k\neq n$}\\
V(\omega_n)\oplus V(\ol{\omega}_n) & &\mbox{ for $D_n$ and $k=n$}\\
\bigoplus_{l=0}^{[k/2]}V(\omega_{k-2l}) & &\mbox{ for $C_n$}
\end{array}\right.
\]

{\it Example 1}. For $\mathfrak g=
\Aod=A_{7}^{(2)\dagger}$,
$$V^4=V^{\otimes 4}/(V\ot V\ot W+V\ot W\ot V+ W\ot V\ot V)
\simeq V(\omega_4)\oplus V(\overline{\omega}_4),$$
\noindent where the highest weight vectors are
$u_{\omega_4}=v_1\wedge v_2\wedge v_3\wedge v_4$ and
$u_{\overline{\omega}_4}=v_1\wedge v_2\wedge v_3\wedge 
v_{-4}$.
$$V^3
=V^{\otimes 3}/(V\otimes W+W\otimes V)\simeq V(\omega_3),$$
where the highest weight vector is $v_1\wedge v_2\wedge
v_3.$ 
Note that
$$
e_3(v_1\wedge v_2\wedge v_{-2}-q v_1\wedge v_3\wedge v_{-3}
+q^2 v_1\wedge v_4\wedge v_{-4})
=-q(1+q^2)v_1\wedge v_3\wedge v_4.
$$

$$V^2
=V^{\otimes 2}/W\simeq V(\omega_2),$$
where the highest weight vector is $v_1\wedge v_2$
and we have
\begin{eqnarray*}
&&v_1\wedge v_{-1}-qv_2\wedge v_{-2}+q^2v_3\wedge v_{-3}
-q^3v_4\wedge v_{-4}\\
&=&-q^6(v_{-1}\wedge v_{1}-q^{-1}v_{-2}\wedge v_{2}+q^{-2}v_{-3}\wedge v_{3}
-q^{-3}v_{-4}\wedge v_{4})
\end{eqnarray*}

{\it Example 2.} For $\mathfrak g= \Ao=A_{5}^{(2)}$,
$$V^2=V^{\otimes 2}/W\simeq V(\omega_2)\oplus V(0),$$
\noindent where the highest weight vectors are $u_{\omega_2}=v_1\wedge v_2$ and
$u_{0}=v_1\wedge v_{-1}-qv_2\wedge v_{-2}+q^2v_3\wedge v_{-3}$.
$$V^3
=V^{\otimes 3}/(V\otimes W+W\otimes V)\simeq V(\omega_3)\oplus V(\omega_1),$$
\noindent where the highest weight vectors are
$u_{\omega_3}=v_1\wedge v_2\wedge v_3$ and
$u_{\omega_1}=v_1\wedge v_{2}\wedge v_{-2}-q
v_1\wedge v_3\wedge v_{-3}.$

{\it Example 3.} For $\mathfrak g= C_n^{(1)}=C_4^{(1)}$,
$$V^2=V^{\otimes 2}/W\simeq V(\omega_2)$$
where the highest weight vector is $v_1\wedge v_2$. It is a good exercise
to check directly from the relations of $W$ that
\begin{eqnarray*}
&&v_{-1}\wedge v_1-q^{-1}v_{-2}\wedge v_2+q^{-2}v_{-3}\wedge v_3-
q^{-3}v_{-4}\wedge v_4\\
&=&v_{1}\wedge v_{-1}-qv_{2}\wedge v_{-2}+q^{2}v_{3}\wedge v_{-3}-
q^{3}v_{4}\wedge v_{-4}\\
&=&0
\end{eqnarray*}
In fact denote the two
vectors in the first two lines by $u_0'$ and $u_0$ respectively. Let $W_1$ 
be the submodule generated by the common relations (\ref{E:rel-1}-\ref{E:rel-2}) plus $v_{-1}\wedge v_1=-
v_{1}\wedge v_{-1}$, and $W_2$ to be the submodule generated by 
the common relations (\ref{E:rel-1}-\ref{E:rel-2}) plus $v_{-n}\wedge v_n=-q^2v_{n}\wedge v_{-n}$. 

It can be seen immediately 
that 
$$u_0'\equiv -q^2u_0\ \ mod\, W_2$$
as in the case of $\Ao$ (cf. \ref{E:wt0}). 
Swapping the role of $1$ and $n$ and straightening backward
we have 
$$u_0'\equiv -q^{-(2n-2)}u_0 \ \ mod\, W_1.$$ 
As long as $q^{2n}\neq 1$ we have
$$u_0=u_0'=0. $$ 

Summarizing previous calculations we obtain the following
result.
\begin{proposition}\label{P:relations-1}
Let $V$ be the $2n$ dimensional vector space generated
by basis vectors $v_j$, $j\in J=\{1, 2, \cdots, n, -n, \cdots, -1\}$. Let $W_1$ be the subspace of $V\otimes V$ generated by relations:
\begin{eqnarray*} 
&&v_i\ot v_i,\quad v_i\ot v_j+qv_j\ot v_i~(i\succ j,i\neq\pm j),\\ 
&&v_{-i}\ot v_i+q^2v_i\ot v_{-i}+q(v_{i+1}\ot v_{-i-1}+v_{-i-1}\ot v_{i+1})~
(1\le i<n),\\
&&v_1\ot v_{-1}+v_{-1}\ot v_1.
\end{eqnarray*}
Let $W_2$ be the subspace of $V\otimes V$ with the
same relations except that $v_1\ot v_{-1}+v_{-1}\ot v_1
$ is replaced by $v_{-n}\ot v_{n}+q^2v_{n}\ot v_{-n}$. Then we
have
\begin{eqnarray*}
\sum_{i=1}^{n}(-q)^{-(i-1)}v_{-i}\ot v_{i}
&\equiv &-q^{2}(\sum_{i=1}^{n}(-q)^{i-1}v_{i}\ot v_{-i})\ \ mod\, W_2\\
&\equiv &-q^{-(2n-2)}(\sum_{i=1}^{n}(-q)^{i-1}v_{i}\ot v_{-i})
\ \ mod\, W_1.
\end{eqnarray*}
Also we have
\begin{eqnarray*}
v_{-i}\ot v_{i}
&\equiv &-v_i\ot v_{-i}+(1-q^2)\sum_{k=1}^{i-1}(-q)^{-k}
v_{i-k}\ot v_{-(i-k)}\ \ mod\, W_1 \\
&\equiv &-q^2v_i\ot v_{-i}+(1-q^2)\sum_{k=1}^{n-k}(-q)^{k}
v_{i+k}\ot v_{-(i+k)}\ \ mod\, W_2.\\
\end{eqnarray*}
\end{proposition}

\section{Crystal structures}

In this section we show that the $\Uqp$-modules $V^k$ constructed in
the last section admit crystal base in the
sense of Kashiwara. Since $V^k$ for $\Ao$ and $\Aod$ are the same as
$\Uqp$-modules (apart from the labeling of the vertices of the
Dynkin diagram), we consider the $\Aod$ case. The crystal structure
of $V^k$ for $\Ao$ is rather difficult to describe. Before going
into details we recall necessary facts about crystal base from \cite{K1}.

\subsection{Basics of crystal base}
Let $\ol{P}=\sum_{j=1}^n\Z\ol{\La}_j$ be the weight lattice and
$\ol{P}^+$ denote the set of dominant integral weights
for $\ol{\geh}$. Let $V$ be a finite dimensional $\Uqp$-module.
Then $V=\bigoplus_{\mu\in\ol{P}}V_\mu$ where $V_\mu=\{u\in V\mid
t_iu=q_i^{\mu(h_i)}u\}$. We shall explain the crystal base for $V$.
Kashiwara defined the operators
$\et{i}$ and $\ft{i}$ acting on $V$ by
\[
\et{i}f_i^{(k)}u=f_i^{(k-1)}u,\ft{i}f_i^{(k)}u=f_i^{(k+1)}u
\mbox{ for }u\in V_{\mu+k\alpha_i}\cap\mbox{Ker }e_i.
\]
Let $A$ be the ring of rational functions regular at
$q=0$. The crystal base for $V$ is a pair $(L,B)$ satisfying the
following properties.
\begin{eqnarray}
\bullet\quad&&\mbox{$L$ is a free $A$-submodule of $V$ such that }
\Q\ot_AL\simeq V. \label{eq:cry1}\\
\bullet\quad&&\mbox{$B$ is a base of the $\Q$-vector space }L/qL.\\
\bullet\quad&&\et{i}L\subset L\mbox{ and }\ft{i}L\subset L\mbox{ for any }i.
\label{eq:cry3}\\
&&\mbox{Hence $\et{i}$ and $\ft{i}$ act on $L/qL$}.
\nonumber\\
\bullet\quad&&\et{i}B\subset B\cup\{0\}\mbox{ and }\ft{i}B\subset B\cup\{0\}.
\label{eq:cry4}\\
\bullet\quad&&L=\bigoplus_{\la\in\ol{P}}L_\la\mbox{ and }
B=\bigsqcup_{\la\in\ol{P}}B_\la,\\
&&\mbox{where }L_\la=L\cap V_\la\mbox{ and }B_\la=B\cap(L_\la/qL_\la).
\nonumber\\
\bullet\quad&&\mbox{For $b,b'\in B$, $b'=\ft{i}b$ if and only if
$b=\et{i}b'$}. \label{eq:cry6}
\end{eqnarray}
$L$ is called the crystal lattice and $B$ the crystal
for $V$. The set $B$ has a graph structure which is called the crystal
graph for $V$.

Let $V(\la)$ be the irreducible $U_q(\ol{\geh})$-module with highest weight
$\la\in\ol{P}^+$ and highest weight vector $u_\la$, and set
\begin{eqnarray*}
L(\la)&=&A\langle\ft{i_1}\cd\ft{i_l}u_\la\mid
l\ge0,1\le i_1,\cd,i_l\le n\rangle,\\
B(\la)&=&L(\la)/qL(\la)\setminus\{0\}.
\end{eqnarray*}
Then up to a trivial isomorphism, $(L(\la),B(\la))$ is the unique crystal
base of $V(\la)$. It is also known that if
$V=V(\la_1)\oplus V(\la_2)$ ($\la_1,\la_2\in\ol{P}^+$) and
$(L(\la_i),B(\la_i))$ is the crystal base for $V(\la_i)$
($i=1,2$), then $(L(\la_1)\oplus L(\la_2),B(\la_1)\sqcup B(\la_2))$
is the crystal base for $V$. For $\la\in\ol{P}^+$, the
explicit realization of the crystal $B(\la)$ for the $U_q(\ol{\geh})$-module
$V(\la)$ is given in \cite{KN}.

Our strategy to show that $V^k$ has a crystal is as follows. Since we already
know the structure of $V^k$ as a $U_q(\ol{\geh})$-module, we take
$(L(\omega_k),B(\omega_k))$ (or $(L(\omega_n)\oplus L(\ol{\omega}_n),
B(\omega_n)\sqcup B(\ol{\omega}_n))$ when $\geh=\Aod,k=n$) as a candidate
for crystal base. Then it is clear that the conditions
(\ref{eq:cry1})-(\ref{eq:cry6}) are valid except (\ref{eq:cry3}),
(\ref{eq:cry4}) and (\ref{eq:cry6}) for $i=0$. Therefore, we are left to
check those conditions for $i=0$. We do it with the help of the explicit
description of the crystal structure given in \cite{KN}.

\subsection{Crystal bases}
Let us set
\[
\begin{array}{clll}
&L^k=L(\omega_n)\oplus L(\ol{\omega}_n),
&B^k=B(\omega_n)\sqcup B(\ol{\omega}_n)&\quad\mbox{ for }\ol{\geh}=D_n,k=n,\\
&L^k=L(\omega_k),&B^k=B(\omega_k)&\quad\mbox{ otherwise.}
\end{array}
\]
We shall show that $\et{0},\ft{0}$ preserve $L^k,B^k\cup\{0\}$, and
their actions on $B^k$ is described explicitly. We recall the
description of $B^k$ in \cite{KN}.
\[
B^k=\left\{b_{i_1,i_2,\cd,i_k}\left|
\begin{array}{l}
i_1,\cd,i_k\in J\mbox{ satisfying conditions}\\
\mbox{(1) and (2) below}
\end{array}\right.\right\}.
\]
\begin{itemize}
\item[(1)] $i_\nu\prec i_{\nu+1}$ for $1\le\nu<k$,\quad if $\ol{\geh}=C_n$,\\
$i_\nu\prec i_{\nu+1}$ or $i_\nu=i_{\nu+1}=0$ for $1\le\nu<k$,\quad
if $\ol{\geh}=B_n$,\\
$i_\nu\prec i_{\nu+1}$ or $i_\nu=-i_{\nu+1}=-n$ for $1\le\nu<k$,\quad
if $\ol{\geh}=D_n$.
\item[(2)] If $i_s=p$ and $i_t=-p$, then $s+(k-t+1)\le p$.
\end{itemize}
Let $v_{i_1,i_2,\cd,i_k}$ denote $v_{i_1}\wedge v_{i_2}\wedge\cd\wedge 
v_{i_k}$ in $V^k$ for $i_1,i_2,\cd,i_k\in J$.
By the uniqueness theorem of crystal base, as a $U_q(\ol{\geh})$-crystal,
we identify $(L^k,B^k)$ with the crystal base of $V^k$ by
$b_{1,2,\cd,n-1,\pm n}=v_{1,2,\cd,n-1,\pm n}\mbox{ mod }qL^n$ for
$\ol{\geh}=D_n,k=n$, and $b_{1,2,\cd,k}=v_{1,2,\cd,k}\mbox{ mod }qL^k$
otherwise.

By the explicit actions of $\et{i},\ft{i}$ ($1\le i\le n$) on $B^k$ given in
\cite{KN} we have the following result.

\begin{lemma} \label{lem:1}
We have
\begin{itemize}
\item[(1)] Except for $\geh=\Aod$ and $k=n$, if
           $b_{1,i_2,\cd,i_k}=\ft{a_1}\cd\ft{a_l}b_{1,2,\cd,k}$, then
           $1\notin\{a_1,\cd,a_l\}$ and furthermore,
           $\ft{a_1}\cd\ft{a_l}b_{2,3,\cd,k,-1}=b_{i_2,i_3,\cd,i_k,-1}$.
\item[(2)] For $\geh=\Aod$ and $k=n$, if
           $b_{1,i_2,\cd,i_n}=\ft{a_1}\cd\ft{a_l}b_{1,2,\cd,n-1,\pm n}$, then
           $1\notin\{a_1,\cd,a_l\}$ and furthermore,
           $\ft{a_1}\cd\ft{a_l}b_{2,3,\cd,n-1,\pm n,-1}
            =b_{i_2,i_3,\cd,i_n,-1}$.
\item[(3)] Except for $\geh=\Aod$ and $k=n$, if
           $b_{i_1,\cd,i_{k-1},-1}=\et{a_1}\cd\et{a_l}b_{-k,\cd,-2,-1}$, then
           $1\notin\{a_1,\cd,a_l\}$ and furthermore,
           $\et{a_1}\cd\et{a_l}b_{1,-k,\cd,-2}=b_{1,i_1,\cd,i_{k-1}}$.
\item[(4)] For $\geh=\Aod$ and $k=n$, if
           $b_{i_1,\cd,i_{k-1},-1}
            =\et{a_1}\cd\et{a_l}b_{\pm n,-(n-1),\cd,-2,-1}$, then
           $1\notin\{a_1,\cd,a_l\}$ and furthermore,
           $$\et{a_1}\cd\et{a_l}b_{1,\pm n,-(n-1),\cd,-2}
            =b_{1,i_1,\cd,i_{k-1}}.$$
\end{itemize}
\end{lemma}

Define the following operators.

\begin{eqnarray*}
F_k&=&\ft{1}^2\cd\ft{k-1}^2\ft{k}\cd\ft{n-1}\ft{n}\ft{n-1}\cd\ft{k+1}\ft{k},
\quad\mbox{ for }\geh=\C,\\
&=&\ft{1}^2\cd\ft{k-1}^2\ft{k}\cd\ft{n-1}\ft{n}^2\ft{n-1}\cd\ft{k+1}\ft{k},
\quad\mbox{ for }\geh=\Ae,\\
&=&\ft{1}^2\cd\ft{k-1}^2\ft{k}\cd\ft{n-2}\ft{n}\ft{n-1}\cd\ft{k+1}\ft{k},
\quad\mbox{ for }\geh=\Aod,k\ne n,\\
F_n^+&=&\ft{1}^2\cd\ft{n-2}^2\ft{n}^2,\quad
F_n^-=\ft{1}^2\cd\ft{n-2}^2\ft{n-1}^2,
\quad\mbox{ for }\geh=\Aod,k=n.\\
\end{eqnarray*}
The operators $E_k$ or $E_n^\pm$ are defined similarly by replacing $\ft{i}$
with $\et{i}$. Using the actions of $f_i$ and $e_i$ on $V^k$ we have the
following lemma.

\begin{lemma} \label{lem:2}
We have
\begin{itemize}
\item[(1)] Except for $\geh=\Aod$ and $k=n$, we have
           $F_kv_{1,2,\cd,k}=v_{2,3,\cd,k,-1}$ and
           $f_0v_{2,3,\cd,k,-1}=v_{1,2,\cd,k}$.
\item[(2)] For $\geh=\Aod$ and $k=n$, we have
           $F_n^\pm v_{1,2,\cd,n-1,\pm n}=v_{2,\cd,n-1,\mp n,-1}$ and
           $f_0v_{2,\cd,n-1,\pm n,-1}=v_{1,2,\cd,n-1,\pm n}$.
\item[(3)] Except for $\geh=\Aod$ and $k=n$, we have
           $E_kv_{-k,\cd,-2,-1}=v_{1,-k,\cd,-3,-2}$ and
           $e_0v_{1,-k,\cd,-3,-2}=v_{-k,\cd,-2,-1}$.
\item[(4)] For $\geh=\Aod$ and $k=n$, we have
 \begin{eqnarray*}
E_n^\pm v_{\pm n,-(n-1),\cd,-2,-1}&=&v_{1,\mp n,-(n-1),\cd,-2}\quad\mbox{and}\\
 e_0v_{1,\pm n,-(n-1),\cd,-2}&=&v_{\pm n,-(n-1),\cd,-2,-1}.
\end{eqnarray*}
\end{itemize}
\end{lemma}

Now we have the following theorem declaring the existence of a crystal base for $V^k$.

\begin{theorem}
We have
\begin{itemize}
\item[(1)] $\et{0}$ and $\ft{0}$ preserve $L^k$.
\item[(2)] $\ft{0}b_{i_1,i_2,\cd,i_k}=b_{1,i_1,\cd,i_{k-1}}$ if $i_k=-1$,
           $=0$ otherwise.
\item[(3)] $\et{0}b_{i_1,i_2,\cd,i_k}=b_{i_2,\cd,i_k,-1}$ if $i_1=1$,
           $=0$ otherwise.
\end{itemize}
\end{theorem}
{\it Proof.} Note that 
(2) and (3) ensures (1).

We first prove (2) except for the case of $\geh=\Aod,k=n$. If $i_k\ne-1$,
$\ft{0}b_{i_1,\cd,i_k}=0$ from the weight consideration. If $i_k=-1$, take
$a_1,\cd,a_l\in\{2,\cd,n\}$ such that $b_{1,i_1,\cd,i_{k-1}}=
\ft{a_1}\cd\ft{a_l}b_{1,2,\cd,k}$. Note that $f_0\ft{i}=\ft{i}f_0$ for
$i\ne0,1$ and $b_{2,\cd,k,-1}=v_{2,\cd,k,-1}\mbox{ mod }qL^k$ from
Lemma \ref{lem:2} (1). Using Lemma \ref{lem:1} (1) and Lemma \ref{lem:2} (1),
we have
\begin{eqnarray*}
\ft{0}b_{i_1,\cd,i_{k-1},-1}&=&f_0\ft{a_1}\cd\ft{a_l}b_{2,\cd,k,-1}
=f_0\ft{a_1}\cd\ft{a_l}v_{2,\cd,k,-1}\\
&=&\ft{a_1}\cd\ft{a_l}f_0v_{2,\cd,k,-1}
=\ft{a_1}\cd\ft{a_l}v_{1,2,\cd,k}=b_{1,i_1,\cd,i_{k-1}}
\end{eqnarray*}
mod $qL^k$.

For the case of $\geh=\Aod,k=n$, the proof is similar. The only difference
is that we have to choose $+$ or $-$ of $b_{1,i_1,\cd,i_{n-1}}
=\ft{a_1}\cd\ft{a_l}b_{1,2,\cd,\pm n}$, since $B^n$ has two connected
components.

For the proof of (3), observe that
$b_{\pm n,-(n-1),\cd,-2,-1}=v_{\pm n,-(n-1),\cd,-2,-1}$
mod $qL^n$ for $\geh=\Aod,k=n$, and
$b_{-k,\cd,-2,-1}=v_{-k,\cd,-2,-1}$ mod $qL^k$ otherwise.
\hspace{\fill} $\Box$

\vskip0.5cm\noindent
{\bf Acknowledgements}.
This work was started while M.O. was visiting Department of Mathematics,
North Carolina State University in March-April, 1998. M.O. thanks the
staffs for their warm hospitality. 
N.J. acknowledges the partial support from  NSA grant 
MDA904-97-1-0062. 
K.C.M. acknowledges the partial support from  NSA grant MDA904-96-1-0013.

\end{document}